\documentclass[a4paper,10pt]{article}
\usepackage[utf8]{inputenc}
\usepackage[english]{babel}
\usepackage{amssymb}
\usepackage{amsmath}
\usepackage{epsfig}
\usepackage{bm}
\usepackage{graphics}
\usepackage{slashbox}
\def\bl{\rule[-1mm]{2.4mm}{2.4mm}}

\def\be{\begin{equation}}
\def\ee{\end{equation}}

\newtheorem{thrm}{\bf Theorem}
\newtheorem{lmm}{\bf Lemma}
\newtheorem{dfn}{\bf Definition}
\newtheorem{rmk}{\bf Remark}

\begin{document}
\title {Water flow under rectangular dam\thanks{Supported by RSCF grant 16-11-10349}}
\author{A.B.~Bogatyr\"ev \and  O.A.~Grigor'ev}
\date{}
\maketitle
\abstract{
A new analytical method for calculation the characteristics of the flow in porous medium bounded by a rectangular polygonal path 
is proposed. The method is based on the usage of high genus Riemann theta functions.}

{\bf MSC2010:} 30C20, 14H42, 30C30

{\bf Keywords:} rectangular polygons, Schwarz-Christoffel integral, theta functions, conformal mapping, Jacobian 

\section{Introduction}
Numerical conformal mappings are commonplace in 2D filtration theory \cite{Birk}. The majority of model problems being considered deal with flows of a complex kinematic structure, with a lot of singularities (e.g., points where velocity reaches zero or pressure tends to infinity). Therefore, the boundaries of hodograph areas are rich in spikes and right angles. A typical example of such an area is a special case of a circular polygon, namely, a rectangular polygon on a polar grid.

A recent trend in filtration theory is to develop numerical techniques which deal specifically with said types of areas \cite{Ber16,BAP}. These is mostly due to computational problems all the general techniques are ridden with, the most important of those being the crowding phenomenon. However, there already exists an approach addressing precisely these problems \cite{AB,BG2}. The cases treated in \cite{Ber16} are essentially special cases of this general approach.

To illustrate the power of this technique, however, it's sufficient to consider an even simpler case - an underground flow in porous media bounded by a watertight surface from below. Despite such a design may seem unrealistic (as a rule, even the simplest models take groundwater afflux into account, optionally along with a watertight structure covering a finite area), our results may have some practical use when rivers with rock formations underneath are dealt with. This toy model also allows a researcher to solve some geometric design problems, which we illustrate in this paper. Other problems, e.g. free boundary problems and finding a depression curve are yet to be addressed. 

\section{Problem setting}
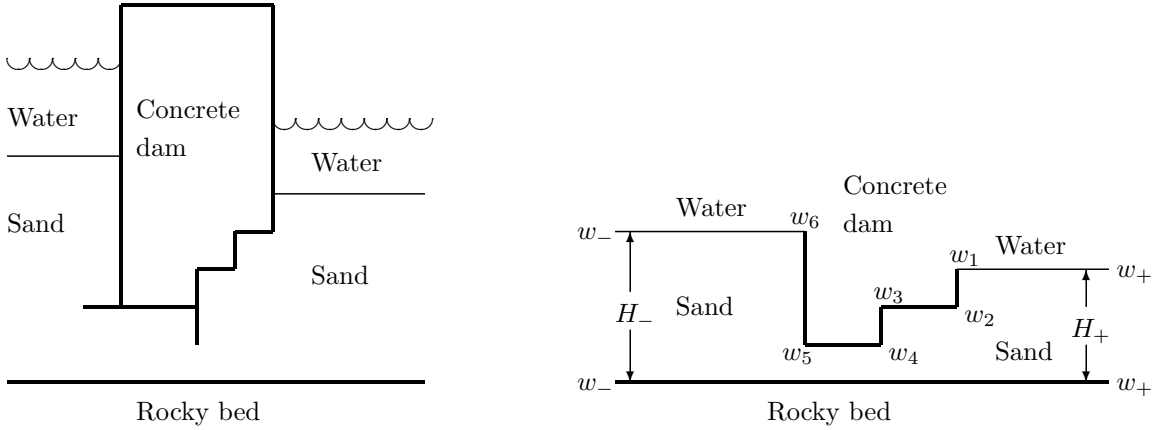
\begin{figure}
\begin{picture}(170,65)
\thicklines
\linethickness{.5mm}
\put(0,5){\line(1,0){55}}
\put(15,15){\line(0,1){40}}
\put(25,10){\line(0,1){10}}
\put(30,20){\line(0,1){5}}
\put(35,25){\line(0,1){30}}
\put(10,15){\line(1,0){15}}
\put(25,20){\line(1,0){5}}
\put(30,25){\line(1,0){5}}
\put(15,55){\line(1,0){20}}
\put(17,40){Concrete}
\put(17,35){dam}
\put(17,0){Rocky bed}
\put(40,18){Sand}
\put(0,25){Sand}
\put(40,33){Water}
\put(0,39){Water}

\thinlines
\linethickness{.2mm}
\put(0,35){\line(1,0){15}}
\put(35,30){\line(1,0){20}}
\multiput(1.5,48)(3,0){5}{\oval(3,3)[b]}
\multiput(36.5,40)(3,0){7}{\oval(3,3)[b]}
\thicklines
\linethickness{.5mm}
\put(80,5){\line(1,0){65}}
\put(105,10){\line(0,1){15}}
\put(115,10){\line(0,1){5}}
\put(125,15){\line(0,1){5}}
\put(105,10){\line(1,0){10}}
\put(115,15){\line(1,0){10}}
\put(110,30){Concrete}
\put(110,25){dam}
\put(100,0){Rocky bed}
\put(130,8){Sand}
\put(88,14){Sand}
\put(130,22){Water}
\put(88,27){Water}

\thinlines
\linethickness{.2mm}
\put(80,25){\line(1,0){25}}
\put(125,20){\line(1,0){20}}
\put(124,21){$w_1$}
\put(126,13){$w_2$}
\put(114,16){$w_3$}
\put(116,8){$w_4$}
\put(102,8){$w_5$}
\put(103,26){$w_6$}
\put(146,4){$w_+$}
\put(146,19){$w_+$}
\put(75,4){$w_-$}
\put(75,24){$w_-$}
\put(80,13){$H_-$}
\put(140,11){$H_+$}
\put(82,17){\vector(0,1){8}}
\put(82,12){\vector(0,-1){7}}
\put(142,15){\vector(0,1){5}}
\put(142,10){\vector(0,-1){5}}

\end{picture}
\caption{Left: Cross section of a rectangular dam; ~~~ 
Right: Porous layer is a domain of the solution }
\label{damSect}
\end{figure}

We consider a dam of rectangular cross section built on a layer of porous medium (sand or clay) with underlying rocky foundation as in Fig. \ref{damSect}. The body of the dam may have intruding riblets. In accordance with the Darcy law the velocity of the stationary underground water flow is proportional to the gradient of pressure: $V=-\varkappa\nabla p$. The continuity condition $\nabla\cdot V=0$ implies that the pressure field is harmonic $\Delta p(w)=0$ in the medium. Equipped with natural boundary conditions, namely the pressure $p(w)$ being locally constant on the bottom of upper/lower basin while its normal derivative being equal to zero on the watertight surfaces of the dam and the rocky foundation, this brings us to the mixed boundary value problem for the pressure $p(w)$.

\hspace{1cm}

Consider an analytic function $f(w):=p+iq$ with $q(w)$ being harmonic conjugate function to $p(w)$ which is locally constant at the surface of the water resistant surfaces. This function conformally maps the porous layer to a rectangle whose aspect ratio is related to the 
drop of the water level behind the dam and to the integral water flow  $Q$ through it. The latter value is one of the main characteristics of the hydraulic structure:

\be
Q:=\int_{-\infty}^0-\varkappa \frac{dp}{dy}~dx=\varkappa\int_{-\infty}^0\frac{dq}{dx}~dx=\varkappa(q(0)-q(-\infty)),
\ee
which is proportional to the height of the rectangle (provided zero is the horizontal endpoint of the upper basin). The width of the rectangle is proportional to the drop of the water level behind the dam 
minus the difference in the levels of the bottom itself. The function $f(w)$ maps stream lines of the flow to the horizontal lines in the rectangle.  

The domain filled by the porous medium will be mapped to the rectangle in two steps. First we 
map it to the upper half plane $\mathbb{H}:=\{x\in\mathbb{C}:~Im~x>0\}$. The second step consists in mapping the half plane to the rectangle 
with four marked points on the boundary (where the type of the boundary conditions breaks) being mapped to the vertexes of the rectangle.
The second map is an elliptic integral, its calculation is a standard procedure -- see e.g. \cite{Akh, LS}

\section{Space of octagonal domains}
First we consider a relatively simple case when the porous medium fills out an
octagonal domain shown in  Fig. \ref{damSect} on the right. After we explain the key 
features of this approach, it will be shown how to increase the complexity of the domain
by adding more rectangular steps and/or vertical or horizontal cuts.

We enumerate the vertexes of the octagon in the counterclockwise direction
starting from the endpoint of the channel pointing in the eastward direction: $+\infty:=w_+, ~w_1, ~w_2, \dots, w_6, w_-:=-\infty$. 
Exactly three of its right angles are intruding, that is equal to $\frac32\pi$ if measured from inside
the domain. We parametrize the space $\cal P$ of polygons  shown in the right Fig. \ref{damSect} with the lengths $H_s$ of their sides, to which we ascribe signs for technical reasons:
\be
\label{HeptDim}
i^sH_s:=w_{s+1}-w_s,\qquad s=1,2,\dots,5.\\
\ee
Obviously, $H_1$ and $H_4$ are negative, the rest of the values $H_s$ are positive. We also consider the positive widths $H_\pm$ of two channels as in
Fig. \ref{damSect} (right). Those seven values satisfy a single linear equation and two more inequalities
\be
H_++H_1-H_3+H_5=H_-;
\qquad\qquad H_++H_1>0; \qquad H_++H_1-H_3>0,
\label{NoDegeneracy} 
\ee
which eliminate the degenerate polygons having e.g. self intersection of the boundary or vanishing isthmi.

\section{Hyperelliptic curves with six real branch points}
The conformal mapping of the upper half plane to any octagon
from the space ${\cal P}$ may be represented by the Schwarz-Christoffel integral: 
\be
\label{SCI}
w(x):=A\int^x\frac{(t-x_2)(t-x_4)(t-x_5)dt}{(t-x_+)(t-x_-)y(t)},
\qquad A>0,
\ee
where the real points $x_*$ are apriori unknown pre-images of the vertexes $w_*$ of the octagon and $y^2(t)=\prod_{j=1}^6(t-x_j)$. 
This integral is naturally defined on a hyperelliptic curve with six real branch points and will be eventually expressed in terms of other function theoretic objects. In this section we list essential facts about these curves. The survey is a short version
of Sects. 2,  4 of \cite{AB} and does not contain any proofs, more details may be found in the textbooks like \cite{FK,GH}. 

\subsection{Algebraic model}
The double cover of the sphere with six real branch points
$x_1<x_2<\dots<x_5<x_6$
is a compact genus two Riemann surface $X$, the coordinates of points in its affine part satisfying
\be
\label{X}
y^2=\prod\limits_{s=1}^6(x-x_s),
\qquad  (x,y)\in\mathbb{C}^2.
\ee
This curve admits a \emph{conformal involution} $J(x,y)=(x,-y)$ with six stationary points $p_s:=(x_s,0)$ and an \emph{anticonformal involution} (or \emph{reflection})
$\bar{J}(x,y)=(\bar{x},\bar{y})$. The stationary points set of the latter has three components known as \emph{real ovals} of the curve. Each real oval is an embedded circle and doubly covers exactly one of the segments  $[x_2,x_3]$, $[x_4,x_5]$, $[x_6,x_1]\ni\infty$ of the extended real line $\hat{\mathbb R}:={\mathbb R}\cup\infty$. The upper half plane $\mathbb{H}:=\{Im~x>0\}$ -- the image of the octagonal domain --
may be lifted to the Riemann surface \eqref{X} in exactly two ways. We choose the one with the positive value of the branch $y(x)$
for $x\in[x_6,\infty)$ on the boundary of the half plane. The pre-images $x_\pm$ of the endpoints of two channels $w_\pm$ are 
lifted respectively to the points $p_\pm$ on the boundary of embedded upper half plane and lie on the third real oval of the surface \eqref{X}.

\subsection{Cycles, differentials, periods}\label{SectBasisCycles}
\begin{figure}
\centerline{
\includegraphics[scale=1.2]{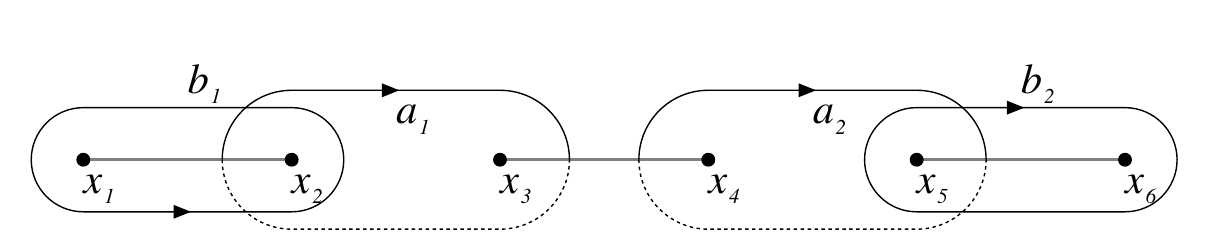}
}
\caption{\small Canonical basis in homologies of the curve $X$}
\label{BasisHomologies}
\end{figure}

We fix a special basis in the 1-homology space of the curve $X$ intrinsically related to the latter. The first and second real ovals give us two 1-cycles,  $a_1$ and $a_2$ respectively. Both cycles are oriented (up to simultaneous change of sign) as the boundary of a pair of pants obtained by removing all real ovals from the surface. Two remaining cycles $b_1$ and $b_2$ are represented by coreal ovals of the curve oriented so that the intersection matrix takes the canonical form --
see Fig. \ref{BasisHomologies}.

The reflection of the surface acts on the introduced basis as follows
\be
\label{abreflect}
\begin{array}{c}
\bar{J}a_s=a_s,
\quad
\bar{J}b_s=-b_s,
\end{array}
\qquad s=1,2.
\ee
Holomorphic differentials on the curve $X$ take the form
\be
\label{diffRep}
du_*=(C_{1*}x+C_{2*})y^{-1}dx,
\ee
with constant values $C_{1*}$, $C_{2*}$.
 The basis of differentials dual to the basis of cycles
\be
\int_{a_s}du_j:=\delta_{sj};
\qquad s,j =1,2,
\ee
determines Riemann period matrix $\Pi$ with the elements
\be
\label{periods}
\Pi_{sj}:=\int_{b_s}du_j;
\qquad s,j =1,2.
\ee
It is a classical result that $\Pi$ is symmetric and has positive definite imaginary part \cite{FK}.

From the symmetry properties (\ref{abreflect}) of the chosen basic cycles it readily follows that:
\begin{itemize}
 \item Normalized differentials are real, i.e. $\bar{J}du_s=\overline{du_s}$, in other words the coefficients $C_*$ in the representation (\ref{diffRep}) are real.
  \item Period matrix is purely imaginary,  therefore we can introduce the symmetric and positive definite real  matrix  $\Omega:=Im(\Pi)$,
\end{itemize}

\subsection{ Jacobian and Abel-Jacobi map }\label{AJimage}
\begin{dfn}
Given a Riemann period matrix $\Pi$, we define the full rank period lattice
$$
L(\Pi)=\Pi\mathbb{Z}^2+ \mathbb{Z}^2 = \int_{H_1(X,\mathbb{Z})}d{\bm u},
\qquad d{\bm u}:=(du_1,du_2)^t,
$$
in $\mathbb{C}^2$ and the 4-torus $Jac(X):=\mathbb{C}^2/L(\Pi)$ known as the Jacobian of the curve $X$.
\end{dfn}
This definition depends on the choice of the basis in the lattice $H_1(X,\mathbb{Z})$, other choices bring us to isomorphic tori.

It is convenient to represent the points  ${\bm u}\in\mathbb{C}^2$ as a theta characteristic $[{\bm\epsilon},{\bm\epsilon}']$, i.e. a 
couple of real 2-vectors (columns)\footnote{Our notation of theta characteristic as two column vectors is not universally accepted: sometimes the transposed matrix is used.}  ${\bm\epsilon}, {\bm\epsilon}'$:
\be
\label{ee'}
{\bm u}=\frac12(\Pi{\bm\epsilon}+{\bm\epsilon}').
\ee
The points of Jacobian  $Jac(X)$ in this notation correspond to
two vectors with real entries modulo 2. Second order points of Jacobian are
represented as $2\times 2$ matrices with binary entries $0,1$.

\begin{dfn}
Abel-Jacobi (briefly: AJ) map 
\be
\label{AJmap}
{\bm u}(p):=\int_{p_1}^p d{\bm u}~~ mod~L(\Pi),
\qquad p_1:=(x_1,0); \quad d{\bm u}:=(du_1,du_2)^t,
\ee
is a correctly defined mapping from the surface $X$ to its Jacobian.
\end{dfn}

From Riemann-Roch formula \cite{FK} it easily follows  that Abel-Jacobi map is a holomorphic embedding of the curve into its Jacobian. In Sect. \ref{theta} we give an explicit equation for the image of the genus two curve in its Jacobian. Meanwhile, let us compute the images of the branching points $p_s=(x_s,0)$, $s=1,\dots,6$, of the curve $X$:

\centerline{
\begin{tabular}{c|c|c}
$p$ & ${\bm u}(p)~ mod~ L(\Pi)$ & $[\epsilon, \epsilon']({\bm u}(p))$\\
\hline
$p_2$ & $\Pi^1/2$ & $\tiny
\left[\begin{array}{c} 10\\00\end{array}\right]$\\
$p_3$&$(\Pi^1+E^1)/2$&$\tiny
\left[\begin{array}{c} 11\\00\end{array}\right]$\\
$p_4$&$(\Pi^2+E^1)/2$&$\tiny
\left[\begin{array}{c} 01\\10\end{array}\right]$\\
$p_5$&$(\Pi^2+E^1+E^2)/2$&$\tiny
\left[\begin{array}{c}01\\11\end{array}\right]$\\
$p_6$&$(E^1+E^2)/2$&$\tiny
\left[\begin{array}{c}01\\01\end{array}\right]$
\label{AJPj}
\end{tabular}}

where $\Pi^s$ and $E^s$ are the $s$-th columns of the period and identity matrix respectively. One can notice that vector ${\bm\epsilon}({\bm u}(p))$ is constant along the real ovals and ${\bm\epsilon}'({\bm u}(p))$ is constant along the coreal ovals.

\subsection{Theta functions on genus two curves}\label{theta}
Here we give a crash course on the theory of Riemann theta functions adapted to genus two surfaces. 
Three problems related to conformal mappings will be effectively solved in terms of Riemann theta functions:
\begin{itemize}
\item Localization of the curve inside its Jacobian;
\item Representation of the 2-sheeted projection of the curve to the sphere;
\item Evaluation of the normalized abelian integrals of the second and the third kinds which appear in the canonical 
decomposition of CS-integral into elementary ones.
\end{itemize}

\begin{dfn}
Let ${\bm u}\in\mathbb{C}^2$ and $\Pi\in\mathbb{C}^{2\times2}$ be a Riemann
matrix, i.e. $\Pi=\Pi^t$ and $Im~\Pi>0$. The theta function of those two arguments is the following Fourier series

$$
\theta({\bm u}, \Pi):=\sum\limits_{{\bm m}\in\mathbb{Z}^2}
\exp(2\pi i{\bm m}^t{\bm u}+\pi i {\bm m}^t\Pi {\bm m}).
$$
Also considered are theta functions with characteristics obtained by a slight modification of the above.
$$
\theta[2{\bm\epsilon}, 2{\bm\epsilon}']({\bm u}, \Pi):=\sum\limits_{{\bm m}\in\mathbb{Z}^2}
\exp(2\pi i({\bm m}+{\bm \epsilon} )^t({\bm u}+{\bm \epsilon}')+\pi i ({\bm m}+{\bm \epsilon})^t\Pi ({\bm m}+{\bm \epsilon}))
$$
$$
=
\exp(i\pi{\bm \epsilon}^t\Pi{\bm \epsilon}+2i\pi{\bm \epsilon}^t({\bm u}+{\bm \epsilon}'))
\theta({\bm u}+\Pi{\bm \epsilon}+{\bm \epsilon}',\Pi),
\qquad {\bm \epsilon},{\bm \epsilon}'\in\mathbb{R}^2.
$$
The matrix argument $\Pi$ of theta function is usually omitted when it is clear which matrix is used.  Omitted vector argument ${\bm u}$ is supposed to be zero and the appropriate function of $\Pi$ is called the theta constant:
$$
\theta[{\bm \epsilon}, {\bm \epsilon}']:=\theta[{\bm \epsilon}, {\bm \epsilon}'](0, \Pi).
$$
\end{dfn}
The convergence of these series follows from $Im~\Pi$ being positive-definite.  The series has high convergence rate with well controlled accuracy \cite{DHB}.
Theta function has the following easily checked quasi-periodicity properties with respect to the lattice $L(\Pi):=\Pi\mathbb{Z}^2+\mathbb{Z}^2$:
\be
\label{quasiperiod}
\theta({\bm u}+\Pi {\bm m}+{\bm m}', \Pi)=\exp(-i\pi {\bm m}^t\Pi {\bm m}-2i\pi {\bm m}^t{\bm u})\theta({\bm u},\Pi),
\qquad {\bm m}, {\bm m}'\in\mathbb{Z}^2.
\ee
Quasi-periodicity of theta with characteristics is easily deduced from the last formula.

\begin{rmk}\label{RemTheta}
(i) Theta function with integer characteristics
$[2{\bm \epsilon}, 2{\bm \epsilon}']$ is either even or odd depending on the parity of the inner product $4{\bm \epsilon}^t\cdot{\bm \epsilon}'$. In particular, all odd theta constants are zeros.\\
(ii) It is convenient to represent integer theta characteristics as the sums of AJ images of the branch points, keeping only the indices of these points: $[sk..l]$ means the sum modulo 2 of the theta characteristics of points
${\bm u}(p_s)$, ${\bm u}(p_k)$, $\dots,{\bm u}(p_l)$, e.g. $[35]$ corresponds to the characteristics
$\tiny \left[\begin{array}{c} 10\\11\end{array}\right]$.
\end{rmk}

\subsubsection{Image of Abel-Jacobi map}
The locus of a genus 2 curve embedded to its Jacobian may be reconstructed by solving a single equation.

\begin{thrm}[Riemann]
A point ${\bm u}$ of Jacobian lies in the image ${\bm u}(X)$ of Abel-Jacobi map
if and only if
\be
\label{JacIm}
\theta[35]({\bm u})=0.
\ee
\end{thrm}
For the localization of the image of AJ embedding for higher genus curves see \cite{AB3}.

\subsubsection{Projection to the sphere}\label{projection}
Any meromorphic function on the curve may be efficiently calculated via the
Riemann theta functions once we know the positions of its zeros and poles (the divisor) \cite{FK,Mum}. 
Take for instance the degree 2 function $x(p)$ on the hyperelliptic curve (\ref{X}). The normalization $x(p_+):=0; x(p_1):=1; x(p_-):=\infty$
leads us to a simple expression for the two sheeted projection to the sphere:
\be
\label{xofP}
x(p):= Const
\frac{\prod_\pm\theta[35]({\bm u\pm \bm u}^+)}{\prod_\pm\theta[35]({\bm u}\pm {\bm u}^-)},
\qquad {\bm u}={\bm u}(p), 
\qquad {\bm u}^\pm:={\bm u}(p_\pm),
\ee
here $\displaystyle{Const:=\frac{\theta^2[35]({\bm u}^-)}{\theta^2[35]({\bm u}^+)}}$.

\subsubsection{Third kind abelian integrals}
On any Riemann surface there exists a unique abelian differential of the third kind $d\eta_{pq}$ with two simple poles at 
the prescribed points $p,~q$ only, residues $+1,-1$ respectively and trivial periods along all $a-$ cycles. There is a closed formula 
for the appropriate abelian integral, due to Riemann: 
\be
\label{kind3}
\eta_{pq}(s):=\int_*^sd\eta_{pq}=
\log\frac{\theta[{\bm\epsilon},{\bm\epsilon}']({\bm u}(s)-{\bm u}(p))}
{\theta[{\bm\epsilon},{\bm\epsilon}']({\bm u}(s)-{\bm u}(q))},
\qquad p,q,s \in X;
\ee
with any odd theta characteristic $[{\bm\epsilon},{\bm\epsilon}']$, say those indicated in the third and the fifth rows of the table in Sect
\ref{AJimage}.

\section{The algorithm of conformal mapping}
The SC integral \eqref{SCI} conformally mapping the upper half-plane  to an octagon from the space $\cal P$ contains the Schwarz-Christoffel differential
of the kind
\be
\label{dw}
dw:=A\frac{(x-x_2)(x-x_4)(x-x_5)dx}{(x-x_+)(x-x_-)y},
\qquad A>0.
\ee
The latter differential  can be canonically decomposed into a sum of $a$-normalized elementary ones:
\be
dw=\frac{H_-}\pi d\eta_{p_-Jp_-}-\frac{H_+}\pi d\eta_{p_+Jp_+} +C_1\cdot du_1+C_2\cdot du_2,
\label{diffdecomp}
\ee
where $d\eta_{pq}$ is an abelian differential of the third kind with poles at $p$ and $q$ only and residues $\pm1$; 
$du_1$ and $du_2$ are holomorphic differentials. All the  differentials here are real, so are all the coefficients in the decomposition. In particular, the residues of $dw$ at its poles are related to the widths of appropriate channels -- this is where the coefficients $H_\pm$ in formula \eqref{diffdecomp} originate from. It remains to get the efficient evaluation of the SC integral itself and all the auxiliary parameters. Using Riemann's  formula (\ref{kind3}), we represent the integral of (\ref{diffdecomp}) in terms of the Jacobian variables:
\be
\label{wab}
w({\bm u})=
\frac{H_-}\pi \log\frac{\theta[3]({\bm u}^-- {\bm u})}{\theta[3]({\bm u}^- + {\bm u})}
-\frac{H_+}\pi \log\frac{\theta[3]({\bm u}^+- {\bm u})}{\theta[3]({\bm u}^+ + {\bm u})}
+ C_1u_1+C_2u_2.
\ee
where  ${\bm u}:=(u_1,u_2)^t={\bm u}(p)$; ${\bm u}^\pm:={\bm u}(p_\pm)=(u_1^\pm,u_2^\pm)$ and the periods matrix for theta is $i\Omega$. 
As such, we found an efficiently computable expression for the SC integral depending on nine real auxiliary parameters contained in $\Omega$, ${\bm u}^\pm$ and $C_1,C_2$.

\subsection{Mapping octagon to the half plane and back}
Here we describe (somewhat sketchy, more details may be found in \cite{AB}) the algorithm of the conformal mapping of a fixed octagon 
from the space ${\cal P}$ to the upper half- plane and back.  Suppose we know the values of all the auxiliary parameters (we postpone their calculation until section \ref{AuxPar}) and a point $w^*$ lies in the octagon normalized by the condition  $w_1=0$,
we consider a system of two complex equations
\be
\label{heptoH}
\begin{array}{r}
w({\bm u}^*) = w^*,\\
\theta[35]({\bm u}^*,i\Omega)=0,
\end{array}
\ee
with respect to the unknown complex 2-vector ${\bm u}^*$. Following arguments in \cite{AB}, we assert that 
system (\ref{heptoH}) has the unique solution ${\bm u}^*$ with theta characteristic from the block
  $\tiny\left[
\begin{array}{cc}
-I~&I\\
-I~&I\\
\end{array}
\right]$, $\qquad I:=(0,1)$.
It is clear (e.g. from the reflection principle for the conformal mappings) that the set of two transcendental equations (\ref{heptoH}) may have many solutions. We use theta characteristic to single out the unique one.

Substituting the solution ${\bm u}^*$ to the right hand side of the expression (\ref{xofP}) we get the evaluation at the point $w^*$ of 
the conformal mapping $x(w)$ of the octagon to the half plane with normalization $w_+$, $w_1$, $w_-\to$  $0,1,\infty$.
In particular, the aspect ratio of the rectangle relating the fall of the water level behind the dam to the total water influx to the upper basin (namely, $\kappa:=\Delta q/\Delta p$) is the ratio of complete elliptic integrals for the curve with the Weierstrass module  
$$1<\lambda:=x(p_6)=\frac{\theta^2[135]({\bm u}^-)\theta^2[356]({\bm u}^+)}
{\theta^2[135]({\bm u}^+)\theta^2[356]({\bm u}^-)}.$$
This value may be eventually expressed via the hypergeometric functions $F(a,b,c|z)$ -- see e.g. \cite{WW} -- as
$$
\kappa=\lambda^{-\frac12}\frac{F(\frac12,\frac12,1|1/\lambda)}{F(\frac12,\frac12,1|1-\lambda)}.
$$

Conversely, given a point $x^*$ in the upper half plane we solve a system of two equations \eqref{xofP} and \eqref{JacIm} with 
respect to a complex 2-vector ${\bm u}^*$ with characteristics from the block we mentioned above. 
Then substitute this solution to the formula (\ref{wab}) for SC integral in the Jacobian variables
to get the image $w^*$ of the point $x^*$ in the octagon.

\subsection{Auxiliary parameters}\label{AuxPar}
Given the dimensions $H_*$ of the octagon, we have to find nine real parameters: the imaginary part $\Omega$ of the $2\times 2$ period matrix, 
the AJ image ${\bm u}^\pm:={\bm u}(p_\pm)$ of two marked points in the Jacobian of the curve, and 2 real coefficients  $C_1$, $C_2$. 
Those parameters give a solution to a system of nine real equations we describe below.

\be
\label{dweq0}
d\theta[35]({\bm u},i\Omega)\wedge
dw({\bm u})=0, 
\qquad {\rm when}~{\bm u}={\bm u}(p_j), \quad j=2,4,5;
\ee
which means that CS differential $dw$ has (double) zeros at three points $p_j$;
\be
\label{ovalu0}
\theta[35]({\bm u}^\pm,i\Omega)=0,
\ee
which means that the points ${\bm u}^\pm$ lie on the AJ image of the curve in the Jacobian and finally
\be\label{Sides}
\begin{array}{l}
-2 H_2=C_1;\\
2 H_4=C_2;\\
-2H_1=2H_-(2u_1^--1)-2H_+(2u_1^+-1)+ C_1\Omega_{11}+C_2\Omega_{12};\\
2H_5=4H_-u_2^--4H_+u_2^+ +C_1\Omega_{12}+C_2\Omega_{22};\\
\end{array}
\ee
Last four equations
specify the side lengths of the polygon and can be deduced by taking the SC integral along the a- and b- cycles 
with subsequent use of the Riemann reciprocity laws \cite{FK,GH} for two latter cases. 
We note that the set of equations \eqref{Sides} immediately gives the unknown coefficients $C_1$, $C_2$ and 
depends linearly on the remaining unknowns. 
\begin{lmm}
Let the side lengths $H_\pm,H_1,H_2\dots,H_5$ satisfy the restrictions \eqref{NoDegeneracy}, then the system of nine equations \eqref{dweq0}, \eqref{ovalu0}, \eqref{Sides} has a unique real solution 
with $0<u_1^+<u_1^-<1/2$ and $\Omega$ in the trihedral cone $0<\Omega_{12}<min(\Omega_{11}, \Omega_{22})$.
\end{lmm}
{\bf Sketch of the proof.} The existence and the uniqueness of the solution to \eqref{dweq0} -- \eqref{Sides} in the specified domain 
follows from the existence and the uniqueness of the conformal mapping of a given octagon to the upper half plane
and explicit description of the image of the periods map \cite{AB}. ~~~\bl

The strategy for the solution of the auxiliary set of equations \eqref{dweq0}-\eqref{Sides}  by Newton method  with parametric continuation is discussed in  \cite{AB}.

\section{Applications}
One can enlarge the stock of polygonal domains at the cost of minor complication of the algorithm.
Suppose three equations (\ref{dweq0}) are not satisfied which means the zeros of the SC differential have moved from the branch points of the curve to the neighbouring real or coreal ovals. Then the composite function  $w({\bm u}(x))$ maps the upper half-plane to an octagon with vertical or horizontal cuts emanating from the intruding right angles. The spaces of the 11-gons of this kind -- with six right angles, three full angles (cuts) and two zero angles at infinity have dimension nine. The corresponding conformal mappings to the half-plane use the moduli space of real genus two curves with five marked points on their (co)real ovals. The parametric representation of the mapping itself is the same as above but it implicitly depends on extra auxiliary parameters. Now we have three more unknown points 
in the Jacobian corresponding to zeros of CS differential and described by six real values and six more real equations. Three of those equations mean the same as \eqref{ovalu0}: additional points live on the AJ images of (co)real ovals of the curve. Another three equations specify the lengths of three cuts. Note that six additional variables do not explicitly participate in the final parametric representation of the conformal mapping, however they influence the values of the involved auxiliary parameters. We hope that an interested reader will easily reconstruct the aforementioned set of 15 equations (many of those being linear) relating 15 unknowns.

Adding new rectangular steps to the domain is more painful for our approach as it rises the genus of the curve bearing the SC integral.
Starting from $g=3$ we have to characterize the period matrices of hyperelliptic curves, which is the simplest version of the notorious 
Schottky problem. For the hyperelliptic case the latter has an effective solution \cite{Mum}:
certain even theta constants vanish at the hyperelliptic locus of the Siegel upper half space.

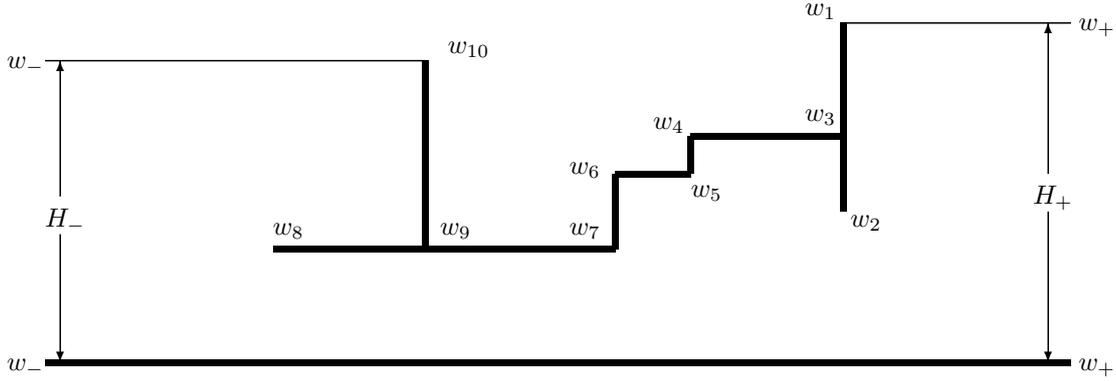
\begin{figure}
	\begin{picture}(170,65)

	\thicklines
	\linethickness{.8mm}
	\put(10,5){\line(1,0){135}}
	\put(60,20){\line(0,1){25}}
	\put(85,20){\line(0,1){10}}
	\put(95,30){\line(0,1){5}}
	\put(40,20){\line(1,0){45}}
	\put(85,30){\line(1,0){10}}
	\put(95,35){\line(1,0){20}}
	\put(115,25){\line(0,1){25}}

	\thinlines
	\linethickness{.2mm}
	\put(10,45){\line(1,0){50}}

	\put(115,50){\line(1,0){30}}
	\put(110,51){$w_1$}
	\put(110,37){$w_3$}
	\put(116,23){$w_2$}
	\put(79,22){$w_7$}
	\put(95,27){$w_5$}
	\put(79,30){$w_6$}
	\put(90,36){$w_4$}
	\put(62,22){$w_9$}
	\put(40,22){$w_8$}
	\put(63,46){$w_{10}$}
	\put(146,4){$w_+$}
	\put(146,49){$w_+$}
	\put(5,4){$w_-$}
	\put(5,44){$w_-$}
	\put(10,23){$H_-$}
	\put(140,26){$H_+$}
	\put(12,27){\vector(0,1){18}}
	\put(12,22){\vector(0,-1){17}}
	\put(142,30){\vector(0,1){20}}
	\put(142,25){\vector(0,-1){20}}

	\end{picture}
	\caption{Example of a ribbed channel corresponding to a curve of genus $g=3$}
	\label{cha1}
\end{figure}

Below we present an example of a construction with two riblets (see fig. \ref{cha1}) along with certain numerical results. First of all, by changing the lengths of the cuts and leaving the rest of the dimensions unchanged (namely, $H_-=13,~H_+=11,~|w_{10}-w_{9}|=10,~|w_{7}-w_{9}|=3,~|w_6-w_7|=|w_5-w_6|=2,~|w_4-w_5|=4,|w_4-w_3|=|w_1-w_3|=2$), we obtain different values for the ratio $\kappa$ of the rectangles' sides we map the porous area to. Those values are summarized in Tab. \ref{Tab1}.

\begin{table}
\centering
	\begin{tabular}{|c|*{5}{c|}}\hline
		\backslashbox{$|w_3-w_2|$}{$|w_9-w_8|$}
		&\makebox[3em]{1}&\makebox[3em]{2}&\makebox[3em]{3}
		&\makebox[3em]{4}\\\hline
		1 & 4.604358822 & 4.866760353 & 5.152737177 & 5.452678199 \\\hline
		2 & 4.647433395 & 4.909830312 & 5.195805421 & 5.495745840 \\\hline
		3 & 4.710043569 & 4.972431551 & 5.258403366 & 5.558342600 \\\hline
		4 & 4.809684090 & 5.072053634 & 5.358018625 & 5.657955418 \\\hline
\end{tabular}
\caption{Aspect ratio $\kappa$ as a function of the lengths of additional edges.}
\label{Tab1}
\end{table}

\begin{figure}
\centering
\includegraphics[keepaspectratio=true, scale=.5]{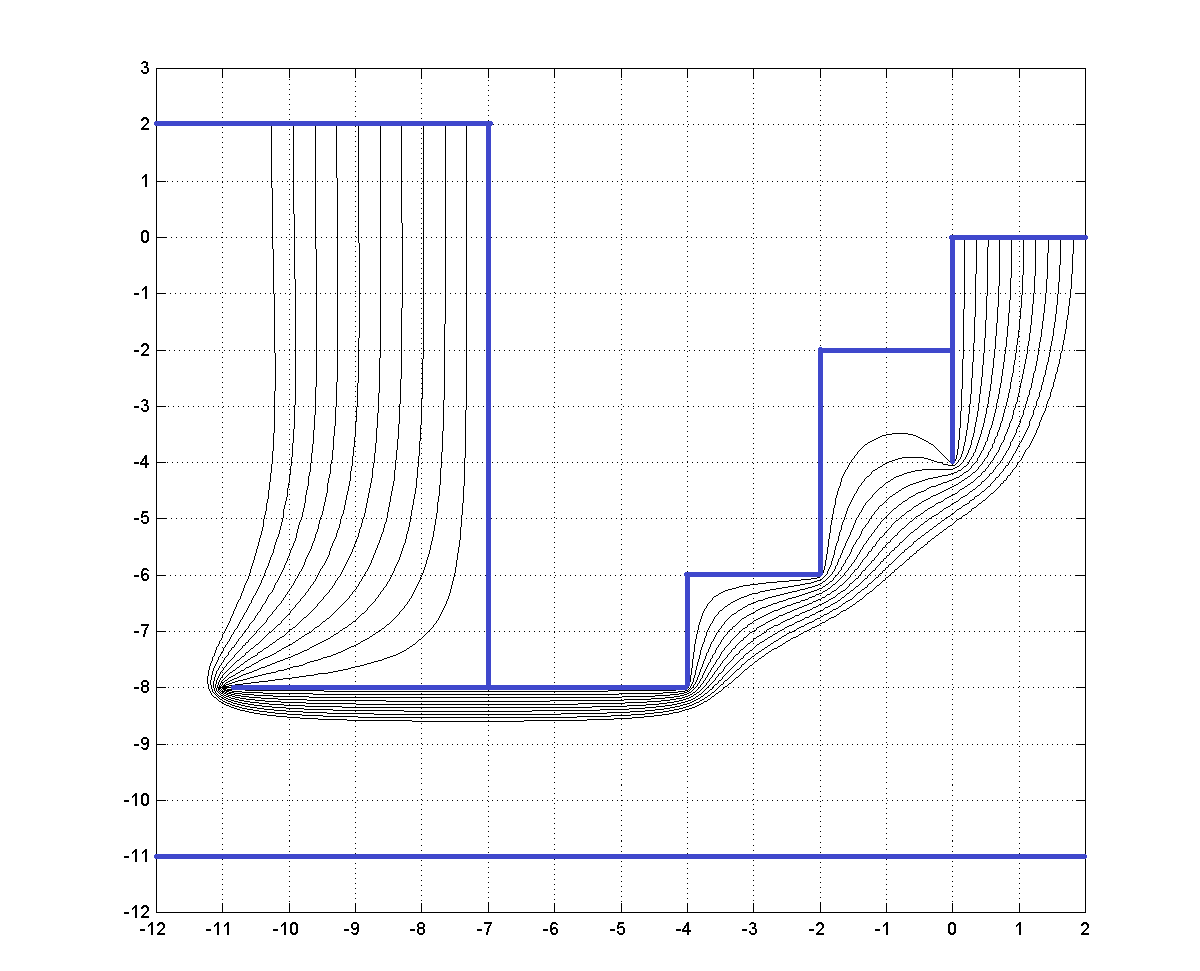}
\caption{Global behaviour of the stream lines under the dam with additional riblets.}
\label{channel1}
\end{figure}

We also chose one in the stock of these polygons ($|w_3-w_2|=2, |w_9-w_8|=4$) to illustrate the global behaviour of stream lines under the dam (Fig. \ref{channel1})

\section{Conclusion}
Today there are many techniques for the realistic numerical simulation of sophisticated 3D filtration problems.
However 2D problems make up an essential part of the latter as main underground stream occurs in the direction across the dam. Coarse characteristic of the flow may be found in the framework of  planar models.
The approach we use reduces finding main characteristics of the underground flow to the solution of a compact set of transcendental and linear equations. This may be applied e.g. to the effective solution of the problems of optimization of the geometry of the dam.

\vspace{5mm}
\parbox{9cm}
{\it
119991 Russia, Moscow GSP-1, ul. Gubkina 8,\\
Institute for Numerical Mathematics,\\
Russian Academy of Sciences\\[3mm]
{\tt ab.bogatyrev@gmail.com, guelpho@mail.ru}}
\end{document}